\documentclass[12pt]{amsart}

\usepackage{amsmath,amsfonts,amsthm,amsopn,cite,mathrsfs}
\usepackage{epsfig,verbatim}
\usepackage{subfigure}

\newcommand{\tfa}{time-frequency analysis}

\newcommand{\stft}{short-time Fourier transform}

\newcommand{\tf}{time-frequency}

\newcommand{\fif}{if and only if}
\newcommand{\tfs}{time-frequency shift}

\newcommand{\modsp}{modulation space}
\newcommand{\psdo}{pseudodifferential operator}

\newcommand{\blt}{Balian--Low theorem}

\newcommand{\knc}{Kohn--Nirenberg correspondence}

\newtheorem{tm}{Theorem}[section]
\newtheorem{lemma}[tm]{Lemma}

\newtheorem{cor}[tm]{Corollary}

\newcommand{\rems}{\noindent\textsl{REMARKS:}}
\newcommand{\rem}{\noindent\textsl{REMARK:}}

 \theoremstyle{definition}
 \newtheorem{definition}{Definition}

\newcommand{\beqa}{\begin{eqnarray*}}
\newcommand{\eeqa}{\end{eqnarray*}}

\newcommand{\field}[1]{\mathbb{#1}}
\newcommand{\bR}{\field{R}}        
\newcommand{\bZ}{\field{Z}}        
\newcommand{\bT}{\field{T}}        %
 \def\cS{\mathcal{S}}
 
 \def\cH{\mathcal{H}}
 \def\cB{\mathcal{B}}
 
 \def\cG{\mathcal{G}}
 \def\cM{\mathcal{M}}

 \def\cA{\mathcal{A}}

 \def\cC{\mathcal{C}}
 
 \def\cO{\mathcal{O}}
 
 \def\cR{\mathcal{R}}

 \def\cV{\mathcal{V}}

\def\vgf{V_gf}

\def\rd{\bR^d}

\def\rdd{{\bR^{2d}}}
\def\zdd{{\bZ^{2d}}}

\def\lrd{L^2(\rd)}

\def\zd{\bZ^d}

\def\mvv{M_v^1}

\def\intrd{\int_{\rd}}
\def\intrdd{\int_{\rdd}}

\def\<{\left<}
\def\>{\right>}

\def\inv{^{-1}}
\def\ud{\, d}

\def\mv1{M_v^1}

\def\Lmpq{L_m^{p,q}}

\def\Mmpq{M_m^{p,q}}
\def\phas{(x,\omega )}

\newcommand{\abs}[1]{\lvert#1\rvert}
\newcommand{\norm}[1]{\lVert#1\rVert}

\newcommand{\weyl}{\sigma ^w}

\newcommand{\fourd}{\bR ^{4d}}
\newcommand{\mif}{M^{\infty,1}}
\newcommand{\phase}{(x,\xi)}

\newcommand{\sjo}{Sj\"ostrand}
\newcommand{\vf}{\varphi}
\newcommand{\vv}{v\circ j \inv}
\newcommand{\lpqm}{\ell ^{p,q}_m}
\renewcommand{\mv}{\cM _v}

\begin{document}
\begin{abstract}
We investigate the properties an exotic  symbol class of \psdo s,
\sjo 's class, with methods   of \tfa\ (phase space
analysis). Compared to the classical treatment, the \tf\ approach
leads to striklingly simple proofs of \sjo 's fundamental results and
to far-reaching generalizations. 
\end{abstract}

\title{Time-Frequency Analysis of Sj\"ostrand's Class}
\author{Karlheinz Gr\"ochenig}
\address{Department of Biomathematics and Biometry, 
 GSF National Research Center for Environment and Health,
 Ingolst\"adter Landstra\ss e 1, D-85764 Neuherberg, Germany}  
\email{karlheinz.groechenig@gsf.de}
\date{}
\maketitle

\section{Introduction}

In 1994/95 \sjo\ introduced a  symbol class for \psdo s that
contains the H\"ormander class $S^0_{0,0}$ and also includes non-smooth
symbols.  He  
proved three fundamental results about the $L^2$- boundedness,
the  algebra property, and the Wiener property. This work had
considerable impact on subsequent work in both hard
analysis~\cite{Bou97,Boul99,her01,Toft01,Toft02a,Toft04} and 
\tfa~\cite{CG03,GH99,GH03}.

\sjo 's  definition goes as follows: Let $g\in \cS (\rdd ) $ be
a $C^\infty $-function with compact support satisfying the property
$\sum _{k\in \zdd } g(t-k) = 1, \forall t\in \rdd $. Then a symbol
$\sigma\in \cS ' (\rdd )$ belongs to  $\mif $, the \sjo\ class,  if 
$$
\intrdd \sup _{k\in \zdd }  |(\sigma \cdot g(.
-k))\,\,\widehat{}\,\, (\zeta )| \, d\zeta  < \infty \, .
$$
The Weyl transfrom of a symbol $\sigma (z,\zeta ) $ is defined as  
\begin{equation}
  \label{eq:f2}
\sigma ^w f(x) =
\intrd \sigma \big(\frac{x+y}{2},\xi \big) e^{2\pi i (x-y)\cdot \xi }
f(y) \, dy d\xi \, .
\end{equation}
 \sjo\ proved the following
fundamental  results about the  Weyl transform of a symbol $\sigma \in
\mif (\rdd ) $~\cite{Sjo94,Sjo95}. 

(a) If $\sigma \in \mif (\rdd ) $, then $\sigma ^w$   
 is a bounded operator on
$\lrd $. 

(b) If $\sigma _1, \sigma _2 \in \mif (\rdd )$ and $\tau ^w = \sigma
_1 ^w \sigma _2^w$,  then  $\tau \in \mif (\rdd )  $; thus  $\mif $ is
a (Banach) algebra of \psdo s. 

(c) If $\sigma \in \mif (\rdd )$ and $\sigma ^w $ is invertible on $\lrd $,
then $(\sigma ^w)\inv = \tau ^w$ for some $\tau \in \mif (\rdd )$. This is
the Wiener property of $\mif $. For the classical symbol classes
results of this type go back to Beals~\cite{beals77}. 

The original  proofs of \sjo\ were carried out in the realm of
classical ``hard'' analysis.
This line of investigation was deepened and extended in subsequent
work by Boulkhemair, Herault, and Toft~\cite{Bou97,Boul99,her01,Toft01,Toft02a,Toft04}.

Later it was  discovered  that \sjo 's class $\mif $ is a special case
of a so-called \modsp . The family of \modsp s had be studied in \tfa\
since the 1980s and later  was also  used  in  the  theory of \psdo s.
The action of \psdo s with classical symbols on \modsp s was
investigated by Tachizawa~\cite{tachizawa94} in 1994; general \modsp s as
symbol classes for \psdo s were introduced in~\cite{GH99}
independently  of \sjo 's work. This line of  investigation and the
emphasis on 
\tf\ techniques was continued in~\cite{CG03,CR02,GH99,GH03,lab01,lab02}.

To make the connection to \tfa , we introduce the operators of
translation and modulation, 
\begin{equation}
  \label{eq:1}
  T_xf(t) = f(t-x) \quad \quad \text{and} \quad \quad M_\omega f(t) =
  e^{2\pi i \omega \cdot t} f(t), \quad \quad t,x,\omega \in \rd  , 
\end{equation}
and note that 
\begin{equation}
  \label{eq:ne1}
(\sigma \cdot g(\cdot
-z))\widehat{}\, (\zeta ) = \intrdd \sigma (t) \bar{g}(t-z) \, e^{-2\pi i
  \zeta \cdot t} \, dt = \langle \sigma , M_\zeta T_z g\rangle \, .
\end{equation}
This is the so-called \stft . It is not only an important and widely
used   \tf\ representation in signal analysis, but an important object
in the mathematical theory of \tfa . A physicist would use a different
terminology for the same object and speak of  position $z$ ,  momentum
$\zeta $,  and
phase space $\rdd $  instead of time and frequency. 

In view of \eqref{eq:ne1} a distribution belongs to \sjo 's class, if
its \stft\ satisfies the condition $\intrdd \sup _{z\in \rdd }
|\langle \sigma , M_\zeta T_x g\rangle | \, d\zeta <\infty $. More
generally, the \modsp s are defined by imposing a weighted $L^p$-norm
on the \stft .  This class of function
spaces  was  
introduced by H.~G.~Feichtinger in
1983~\cite{feichtinger-modulation1} and \cite{feichtinger81,fei83} and has
been studied extensively.  The \modsp s have turned 
out to be the appropriate function and distribution spaces for many 
problems in \tfa .

The objective of this paper is to give the ``natural'' proofs of \sjo
's results. The  definition of \sjo 's by means of the \stft\ \eqref{eq:ne1}
suggests that the mathematics of translation and modulation operators,
in other words, \tfa , should enter in the proofs. 
Although ``natural'' is a debatable notion in mathematics,
we argue that methods of \tfa\ should simplify the original proofs and
shed new light on \sjo 's results. 
 Currently, several
different proofs exist for  the boundedness and the algebra property,
both in the context of ``hard analysis'' and of \tfa . However,  for the  
Wiener property only \sjo 's original  ``hard analysis'' proof was
known, and it was an open problem to find an alternative proof. 

In the following, we will not only give conceptually new and technically
simple proofs of \sjo's fundamental results,   but we will also obtain
new insights.

  Firstly, \tf\ methods provide  detailed information
on which class of function spaces Weyl transforms with symbols in
$\mif $ act boundedly.

  Secondly,  the \tf\ methods   suggest the
appropriate and maximal generalization of \sjo 's  results (to
weighted \modsp s). Although we restrict our attention to Weyl
transforms and \modsp s on $\rd $, all concepts can be defined on
arbitrary locally compact abelian groups. One may conjecture that \sjo
's results hold for (pseudodifferential) operators on $L^2$ of
locally compact abelian groups as well. In that case \tf\ methods hold
more promise than real analysis methods. 

Thirdly,  we show that Weyl transforms with symbols in \sjo 's class
are  almost diagonalized by  Gabor frames.  This may not be
surprising, because it is well-known that   \psdo s with classical
symbols  are almost diagonalized by wavelet bases and local Fourier
bases~\cite{meyer-2,tachizawa-rochberg98}. What is 
remarkable is that the almost diagonalization property with respect to
Gabor frames  is a characterization of \sjo  's class. 

Finally, the new proof of the Wiener property highlights the
interaction with recent  Banach algebra techniques, in particular the
functional calculus in certain 
matrix algebras.

The remainder of the  paper is divided into three parts. In Section~2
 we introduce the basic definitions and results from \tfa . This
area  has now   reached a  level of sophistication that makes it
possible to approach a subject that is usually the domain of ``hard
analysis''. 
In Section~3   we prove the  almost diagonalization property of \psdo s with
symbols in \sjo 's class. In Section~4  we prove \sjo 's results and
their generalization. With
the background in \tfa\ this part becomes very short. 
We conclude with some remarks and problems.

\emph{Acknowledgment.} The author would like to thank Thomas Strohmer
for stimulating  discussions and providing an early version of the
manuscript~\cite{strohmer04},  and furthermore   the Institute for
Mathematical Sciences, National 
University of Singapore, for its hospitality.  This work was finished
at IMS   while the author was visiting   in 2004.

\section{Tools from Time-Frequency  Analysis}

We prepare the tools from \tfa . Most of these are standard and
discussed at length in the text books~\cite{folland89,book}, but the
orginal ideas go back  much further.

\subsection{Time-Frequency Representations}

We combine time $x\in \rd $ and frequency $\xi\in \rd $ into a
single point $z= (x,\xi ) $ in the ``\tf `` plane $\rdd $. Likewise we
combine the operators of translation and modulation to a \emph{\tfs }
and write 
$$
\pi (z ) f(t) = M_\xi T_x f(t) = e^{2\pi i \xi \cdot t}   f(t-x)  
$$
The   \emph{short-time Fourier transform (STFT)} of function/distribution
  $f$ on $\rd $  with   respect to window $g$ is defined by 
  \begin{eqnarray*}
\vgf \phas & = & \intrd f(t) \bar{g}(t-x) e^{-2\pi i t\cdot x} \, dt \\
&=&  \langle f, M_\xi T_x g\rangle = \langle f, \pi (z) g\rangle \, .    
  \end{eqnarray*}
The \stft\ of a symbol $\sigma \phase , \phase \in \rdd ,$ is a
function on $\fourd $ and will be denoted by $\cV _\Phi \sigma
(z,\zeta )$ for $z,\zeta \in \rdd $ in order to distinguish it from
the STFT of a function on $\rd $. 

Usually we fix $g$ in a space of test functions, e.g., $g \in \cS (\rd
)$,  and interpret $f \to V_gf $ as a linear mapping  and $\vgf \phase $
as the \tf\ content of $f$ near the point $\phase $ in the \tf\
plane. 

Similarly, the  (cross-) \emph{Wigner distribution} of $f,g \in \lrd $
is defined as 
\begin{equation*}
  \label{eq:1a}
    W(f,g)(x,\xi)=\int f(x+\frac{t}2)\overline{g(x-\frac{t}2)} e^{-2\pi
    i\xi t}\,dt. 
\end{equation*}
Writing $\check{g}(t) = g(-t)$ for the inversion, we find that the Wigner
distribution is just a \stft\ in disguise:
$$
W(f,g)\phase = 2^d \, e^{4\pi i x\cdot \xi } \, V_{\check{g}}f (2x,2\xi ) \, .
$$
We will need  a well-known  intertwining property
of Wigner distribution, which expresses the Wigner 
distribution of a \tfs\  as  a \tfs , see \cite[p.~57]{folland89}
and~\cite[Prop.~4.3.2]{book}.

\begin{lemma} \label{intertw}
Let $z=(z_1,z_2), w=(w_1,w_2) \in \rdd $ and $f,g\in \lrd $. Then 
\begin{eqnarray*}
& & W(\pi (w) f, \pi (z) g)(x,\xi)\\
&=& e^{\pi i(z_1+w_1)\cdot(z_2-w_2)}\,  e^{2\pi ix\cdot(w_2-z_2 )}
\, e^{2\pi i\xi \cdot(-w_1+z_1)}\, \cdot W(f,g)(x-\frac {w_1+z_1}{2},\xi-\frac
{w_2+z_2 }{2}) \, .
\end{eqnarray*}
In short, with the notation $j(z) = j(z_1,z_2) = (z_2, -z_1)$ we have 
\begin{equation}
  \label{eq:1b}
  W(\pi (w) f, \pi (z) g) = c M_{j(w-z)} T_{\frac{w+z}{2}} W(f,g) \, ,
  \end{equation}
and  the  phase factor $c= e^{\pi i(z_1+w_1)\cdot(z_2-w_2)}$ is  of modulus
$1$.   
\end{lemma}

\subsection{Weyl Transforms}

Using the Wigner distribution, we can recast the definition of the
Weyl transform as follows:
  \begin{equation}
    \label{weyl}
\langle \sigma ^w f, g\rangle =     \langle \sigma , W(g,f)\rangle
\quad \quad f,g \in \cS (\rd )
  \end{equation}
In the context of \tfa\ this is the appropriate definition of the Weyl
transform,  and
we will never use the explicit formula~\eqref{eq:f2}. Whereas the integral
in~\eqref{eq:f2} is defined only for a restricted class of symbols
($\sigma $ should be locally integrable at least), the \tf\ definition
of $\sigma ^w$ makes sense for arbitrary $\sigma \in \cS ' (\rdd
)$.  In addition, if  $T : \cS (\rd ) \to
\cS ' (\rd )$ is continuous, then    the Schwartz
kernel theorem implies that   there exists a $\sigma \in \cS '
(\rdd )$ such that $\langle Tf,g\rangle = (\sigma ^w f, g\rangle $ for
all $f, g \in \cS (\rd )$. Thus,  in a distributional sense,  every
reasonable operator possesses a Weyl symbol. 

The composition of Weyl transforms  defines  bilinear form on 
symbols (\emph{twisted product})
\begin{equation*}
  \label{eq:3}
  \sigma ^w \tau ^w = (\sigma \, \sharp \,   \tau ) ^w
\end{equation*}
Again, there is a (complicated) explicit formula for the twisted
product~\cite{folland89,hormander3}, but 
it is unnecessary for our purpose.

\subsection{Weight Functions}
We use two classes of weight functions. By  $v$ we always  denote  a
non-negative function   on
$\rdd $ with the following properties: \\
(i) $v$ is continuous, $v(0) = 1$,   and  $v$  is  even in each
coordinate $ v(\pm z_1, \pm z_2, \dots , \pm z_{2d}) = v(z_1, \dots , z_{2d})$, \\
(ii) $v$ is  submultiplicative, i.e.,
$v(w+z) \leq v(w) v(z), w,z\in \rdd $, \\
(iii) $v$ satisfies the GRS-condition
(Gelfand-Raikov-Shilov~\cite{gelfandraikov})
\begin{equation}
  \label{eq:2}
  \lim _{n\to \infty } v(nz) ^{1/n} = 1, \quad \quad \forall z \in \rdd
  \, .
\end{equation}
We call a weight satisfying properties (i) --- (iii)
\emph{admissible}. 
Every weight of the form $v(z) =   e^{a|z|^b} (1+|z|)^s \log ^r(e+|z|)
$ for parameters $a,r,s\geq 0$, $0\leq b < 1$ is admissible,  whereas
the exponential weight $v(z) = e^{a|z|}, a>0,$ is
not, because it violates \eqref{eq:2}.

Associated to an admissible weight $v$ ,  we define the class of
so-called \emph{$v$-moderate} weights by 
\begin{equation}
  \label{eq:3b}
  \cM _v = \{ m \geq 0: \sup _{w \in \rdd } \frac{m(w+z)}{m(w)} \leq C
  v(z) , \forall z\in \rdd\} \, .
\end{equation}
Compare also~\cite[Ch.~18.5]{hormander3}. This definition implies that
the weighted mixed-norm $\ell ^p$-space 
$\lpqm $ is  invariant under 
translation whenever  $m\in \mv $. Precisely, set  
$$\|c\|_{\lpqm } =
\Big( \sum _{l\in \zd } \Big(\sum _{k\in  \zd  }
|c_{kl} |^p\, m(\alpha  k, \beta l )^p \Big)^{q/p}\Big)^{1/q}\, ,$$
  and
$(T_{(r,s)} c)_{(k,l)} = c_{(k-r,l-s)}, k,l,r,s \in \zd,$ then
$\|T_{(r,s)} c\|_{\lpqm } \leq C v(\alpha r, \beta s )
\|c\|_{\lpqm }$. Consequently, Young's theorem for convolution implies
that  $\ell ^1_v  \ast \lpqm  \subseteq \lpqm$. 

\subsection{Modulations Spaces and Symbol Classes}

Let $\vf (t) = e^{-\pi t^2}$ be the Gaussian on $\rd $, then we define
a norm on $f$ by imposing a norm on the \stft\ of $f$ as follows:
 \begin{eqnarray*}
     \lefteqn{ \norm{f}_{\Mmpq}=\norm{V_{\vf} f}_{\Lmpq}}\\
& =& \biggl( \intrd \biggl(
    \intrd \abs{\vgf (x,\xi)}^p\, 
m(x,\xi)^p\ud  x\biggr)^{q/p} \ud \xi \biggr)^{1/q}
 \end{eqnarray*}
If  $1\leq p,q< \infty $ and $m\in \mv $, we define $\Mmpq (\rd )  $ as the
 completion of the subspace $\cH _0 = \mathrm{span}\, \{ \pi (z) \vf :
 z \in \rdd \}$ with respect to this norm, if $p=\infty$ or $q=
 \infty $, we use a weak-$^*$ completion. For $p=q$  we write $M^p_m$
 for $M^{p,p}_m$, for $m\equiv 1$, we write $M^{p,q}$ instead of
 $\Mmpq$. For the theory of \modsp s and some applications we refer
 the reader to ~\cite{fei83} and~\cite[Ch.~11-13]{book}

\rems\ 1. The cautionary definition is necessary only for 
weights of superpolynomial growth.  If $m(z) = \cO (|z|^N)$ for some
$N>0$, then $\Mmpq $ is in fact the   subspace of
tempered distributions $f\in \cS ' (\rd )$ for which $\|f\|_{\Mmpq} $
is finite. If $m\geq 1$ and $1\leq p,q \leq 2$, then $\Mmpq $ is a
subspace of $\lrd $. 
However, if  $v(z) = e^{a|z|^b}, b<1$, then $\mvv
\subseteq \cS (\rd )$ and  $\cS'(\rd ) \subseteq M^\infty _{1/v}$, and
we would have to use ultradistributions in the sense of
Bj\"ork~\cite{bjork66} to define $\Mmpq $ as 
a subspace of ``something''. 

2. \emph{Equivalent norms:} Assume that $m \in \mv $ and that $g \in
\mvv $, then 
\begin{equation}
  \label{eq:4}
  \|V_gf \|_{\Lmpq } \asymp \|f\|_{\Mmpq} \, .
\end{equation}
Therefore we  can use arbitrary windows in $\mvv $ in place of the Gaussian
to measure the norm of  $\Mmpq $~\cite[Ch.~11]{book}. In the following we will use this
norm equivalence frequently without mentioning. 

3. The class of \modsp s contains a number of classical function
spaces~\cite[Prop.11.3.1]{book}, in particular $M^{2} = L^2 $; if
$m(x,\xi ) = (1+|\xi 
|^2)^{s/2},s \in \bR $, then $M^2_m = H^2$, the Bessel potential
space; likewise, the Shubin class $Q_s$ can be identified as a
\modsp~\cite{shubin91,BCG04}; 
and even $\cS $ can be represented as an intersection of \modsp s. 

4. If $m \in \mv $, the following embeddings hold for $1\leq p,q
\leq \infty $:
$$
\mvv \hookrightarrow \Mmpq  \hookrightarrow M^\infty _{1/v} \, ,
$$
and $\mvv $ is dense in $\Mmpq $ for $p,q < \infty $, and weak-$^*$
dense otherwise. 

5. The original \sjo\ class is $\mif (\rdd )$~\cite{Sjo94,Sjo95}. We will use the
weighted class $\mif _v $ as a symbol class for \psdo s in our
investigation. For explicitness, we recall the norm of $\sigma \in
\mif _v$:
\begin{equation}
  \label{eq:5}
  \|\sigma \|_{\mif _v } = \intrdd \sup _{z\in \rdd } |\cV _\Phi
  \sigma (z,\zeta )| \, v(\zeta ) \,  d\zeta \, .
\end{equation}
 
In the last few years \modsp s have been used implicitly and explicity
as symbol classes by many authors,
see~\cite{Bou97,Boul99,CG03,CR02,GH99,GH03,GH04,heil-ramanathan97,her01,Toft01,Toft02a,Toft04,lab01,lab02,PT04,tachizawa-rochberg98,tachizawa94}
for a sample of work.




\subsection{Gabor Frames}

Fix a function $g\in \lrd $ and a lattice $\Lambda \subseteq \rdd
$. Usually we take $\Lambda = \alpha \zd \times \beta \zd $ or
$\Lambda = \alpha \zdd $ for some $\alpha ,\beta > 0$. Let $
\cG (g, \Lambda )= \{ \pi (\lambda ) g: \lambda \in \Lambda \}$ be the
orbit of $g$ under $\pi (\Lambda )$.
Associated to $\cG (g,\Lambda) $ we define two operators; first the 
coefficient operator $C_g$ which maps functions to sequences on
$\Lambda $ and is defined by 
\begin{equation}
  \label{eq:13}
  C_gf (\lambda ) = \langle f, \pi (\lambda ) g\rangle \, , \quad \quad
  \lambda \in \Lambda \, ,
\end{equation}
and then the  Gabor frame operator $S = S_{g,\Lambda }$
\begin{equation}
  \label{eq:14}
  Sf = \sum _{\lambda \in \Lambda } \langle f, \pi (\lambda )g\rangle
  \, \pi (\lambda ) g = C_g^* C_g f \, .
\end{equation}
\begin{definition}
  The  set $\cG (g, \Lambda )$ is called a Gabor frame
  (Weyl-Heisenberg frame) for $\lrd $, if $S_{g,\Lambda }$ is bounded and
  invertible on $\lrd $. Equivalently, $C_g$ is bounded from $\lrd $
  into $\ell ^2(\Lambda )$ with \emph{closed range}, i.e. $\|f\|_2
  \asymp \|C_g f \|_2$. 
\end{definition}

If $\cG (g,\Lambda )$ is a frame, then the function  $\gamma = S\inv g\in
\lrd $ is well defined and is called the  (canonical) dual
window. Likewise the ``dual tight frame window'' $\tilde{\gamma } =
S^{-1/2}g$ is in $\lrd $.  Using different  factorizations of the identity
  and the commutativity $S_{g,\Lambda }\pi (\lambda ) =
\pi (\lambda ) S_{g,\Lambda }$ for all $\lambda \in \Lambda $, we
obtain the following series expansions (Gabor expansions) for $f\in \lrd $:
\begin{eqnarray}
  \label{eq:15}
  f&=&  S\inv S = \sum _{\lambda \in \Lambda } \langle f, \pi (\lambda
  )g\rangle \,  \pi (\lambda ) \gamma \\
 &=& S S\inv f = \sum _{\lambda \in \Lambda } \langle f, \pi (\lambda
  )\gamma\rangle \, \pi (\lambda ) g \, . \label{eq:16} \\
&=&  S^{-1/2} S S^{-1/2} f =   \sum _{\lambda \in \Lambda } \langle f,
\pi (\lambda   )S^{-1/2} g\rangle \, \pi (\lambda ) S^{-1/2} g \,
. \label{eq:16b} 
\end{eqnarray}

The so-called ``tight Gabor frame
expansion''~\eqref{eq:16b} is particularly useful and convenient,
because it uses only one window $S^{-1/2}g$ and  behaves like an orthonormal
expansion (with the exception that the coefficients are not unique).  

The existence and construction of Gabor frames for  separable
lattices $\Lambda = \alpha \zd \times \beta \zd $) is  well understood
(see~\cite{daubechies90,book,janssen95,walnut93}) and
 we may  take the existence of
Gabor frames with suitable $g$ for granted.

 The expansions~\eqref{eq:15} -- \eqref{eq:16b} converge unconditionally  in $\lrd $, but
for ``nice'' windows the convergence can be extended to other function
spaces. 

The following theorem summarizes the main properties of Gabor
expansions and the characterization of \tf\ behavior by means of Gabor
frames~\cite{fg97jfa,GL03}

\begin{tm} \label{tfmain}
  Let $v$ be an admissible weight function (in particular $v$
  satisfies the
  GRS-condition~\eqref{eq:2}). Assume that  $\cG (g, \alpha \zd \times
  \beta \zd )$   is a Gabor frame for $\lrd $ and that $g\in \mvv $. Then

(i) The dual window  $\gamma = S\inv g$ and $S^{-1/2}g$ are  also in $  \mvv $.

(ii) If $f\in \Mmpq $, then the Gabor expansions~\eqref{eq:15}
--~\eqref{eq:16b} converge 
  unconditionally in $\Mmpq$ for $1\leq p,q<\infty $ and all $m\in \cM
  _v$, and weak-$^*$
  unconditionally if $p=\infty $ or $q=\infty $. 

(iii) The following norms are equivalent on $\Mmpq$:
\begin{equation}
  \label{eq:17}
  \|f\|_{\Mmpq} \asymp \|C_g f\|_{\lpqm  }\asymp \|C_\gamma f\|_{\lpqm
  } \, .
\end{equation}
\end{tm}

\rem\ When  $g\in M^1\supseteq \mvv $, then $\cG (g,\Lambda )$ is necessarily
overcomplete by the  \blt ~\cite{BHW95}. Although the coefficients $\langle
f, \pi (\lambda )g\rangle $ and $\langle f, \pi (\lambda ) \gamma
\rangle$ are not unique, they are the most convenient ones for \tf\
estimates. 

\section{Almost Diagonalization of Pseudodifferential Operators}

 The tools of the previous section have been developed mainly for
 applications in signal analysis, but in view of the definition of  the Weyl transform
~\eqref{eq:f2} and of \sjo 's class~\eqref{eq:5}, we can taylor these
methods to  the investigation of \psdo s. It is now ``natural'' to
study $\weyl $  on \tfs s  of a fixed  function (``atom'')  and then  study the
matrix of $\sigma ^w $ with respect to a Gabor frame. This idea  is related to the 
confinement characterization of $\mif$~\cite{Sjo95}, but is  conceptually much
simpler.

\subsection{Almost Diagonalization}

We first establish a simple, but crucial relation between the action
of $\sigma ^w$ on \tfs s and the \stft\ of $\sigma $. Recall that 
$$
j(z_1,z_2) = (z_2, -z_1)  \quad \quad \mathrm{for } 
\,\, z= (z_1, z_2 ) \in \rdd \, . 
$$

\begin{lemma} \label{trans}
  Fix a window  $g\in \mvv $ and  $\Phi = W(g
  , g )$.  Then, for $\sigma \in \mif _{\vv } $,  
  \begin{equation}
    \label{eq:4a}
    \big|\langle \sigma ^w \pi (z)\vf , \pi (w) \vf \rangle \big| =  \Big|V_\Phi
    \sigma \Big(\frac{w+z}{2}, j(w-z)\Big)\Big| = \big| \cV_\Phi     \sigma (u,v)| \, , 
  \end{equation}
and 
\begin{equation}
  \label{eq:5a}
  |\cV_\Phi \sigma (u,v)| = \Big|\langle \sigma ^w  \pi (u-\frac{1}{2}j\inv
  (v)g, \pi (u+\frac{1}{2}j\inv   (v))g \big\rangle \Big| \, 
\end{equation}
for $u,v,w,z\in \rdd $. 
\end{lemma}

\begin{proof}
Note that \eqref{eq:4a} and \eqref{eq:5a} are well-defined, because the
assumption $g\in \mvv $ implies that $\Phi = W(g,g) \in M^1_{1\otimes
  (\vv )} (\rdd )$~\cite[Prop.~2.5]{CG03}, and so the \stft\ $V_\Phi
\sigma $ makes 
sense for $\sigma \in \mif _{\vv }$. 

We use  the \tf\ definition of
  the Weyl transform~\eqref{eq:f2} and the intertwining
  property  Lemma~\ref{intertw}, then 
  \begin{eqnarray}
    \langle \sigma ^w \pi (z)g, \pi (w)g\rangle _{\rd} &=& \langle \sigma,
    W(\pi  (w)g, \pi (z)g) \rangle _{\rdd } \notag \\
&=& \langle \sigma, c M_{j(w-z)} T_{\frac{w+z}{2}} W(g,g)\rangle
\\
&=& \bar{c}\, \cV _{W(g,g)}\sigma \big(\frac{w+z}{2}, j(w-z) \big) \, , \label{eq:17b}
  \end{eqnarray}
where $c$ is a phase factor of modulus one. 

To obtain \eqref{eq:5a}, we set $u= \frac{w+z}{2}$ and $v=
j(w-z)$. Then   $w= u+\frac{1}{2}j\inv (v)$ and $  z=
u-\frac{1}{2}j\inv (v) $, and  reading   formula~\eqref{eq:17b}
backwards yields \eqref{eq:5a}. 
\end{proof}

The next result on almost diagonalization  is  crucial  and 
all properties of the Sj\"ostrand class will follow easily.  

\begin{tm}[Almost Diagonalization] \label{almost}
  Fix a non-zero  $g \in M^1_{v}$ and assume that $\cG (g,
  \Lambda )$ is a Gabor frame for $\lrd $. Then the following
  properties are equivalent.

(i) $\sigma \in \mif _{v\circ j\inv}(\rdd )$.
 
(ii) $\sigma \in \cS ' (\rdd )$ and  there exists a function $H\in L^1_v (\rdd )$ such that
\begin{equation}
  \label{eq:6}
  |\langle \sigma ^w \pi (z) g, \pi (w) g \rangle | \leq H(w-z) \quad
  \forall w,z \in \rdd \, .
\end{equation}

(iii) $\sigma \in \cS ' (\rdd )$ and  there exists  a sequence  $h\in
\ell ^1_v (\Lambda )$ such that 
\begin{equation}
  \label{eq:6b}
  |\langle \sigma ^w \pi (\mu) g, \pi (\lambda) g \rangle | \leq
  h(\lambda -\mu )  \quad
  \forall \lambda, \mu  \in \Lambda \, .
\end{equation}
\end{tm}

\begin{proof}
We first prove  the equivalence (i) $\, \Longleftrightarrow \, $
(ii) by means of  Lemma~\ref{trans}. 

(i) $\, \Longrightarrow \, $ (ii) $\quad $ Assume that $\sigma \in
\mif _{v\circ j\inv }   $ and set 
$$
H_0(v) = \sup _{u\in \rdd } |\cV _\Phi \sigma (u,v)|$$
By  definition of $\mif _{\vv} $ we have  $H_0 \in L^1_{\vv } (\rdd
)$, so  Lemma~\ref{trans} implies  that
\begin{eqnarray}
    |\langle \sigma ^w \pi (z)\vf , \pi (w) \vf \rangle | &=&  \Big|\cV _\Phi
    \sigma \big(\frac{w+z}{2}, j(w-z)\big)\Big| \notag \\
&\leq & \sup _{u\in \rdd } |\cV _\Phi \sigma (u,j(w-z))| \label{aha} \\
&=& H_0(j(w-z)) \, .  \notag
\end{eqnarray}
Since $\|H_0 \circ j\|_{L^1_v} = \|H_0\|_{L^1_{\vv }}<\infty $, we can
  take   $H=
H_0 \circ j\inv \in L^1_{v}(\rdd ) $ as the dominating function
in~\eqref{eq:6}.  

(ii) $\, \Longrightarrow \, $ (i) $\quad $ Conversely, assume that
$\sigma \in \cS ' (\rdd )$ and that 
$\sigma ^w $ is almost diagonalized by the \tfs s $\pi (z)$ with
dominating function $H\in L^1_v (\rdd )$ as in~\eqref{eq:6}.
Using the transition formula~\eqref{eq:5}, we find that  
\begin{eqnarray*}
  |\cV _\Phi \sigma (u,v)| &=& \Big| \Big\langle \sigma ^w  \pi
  (u-\frac{1}{2}j\inv 
  (v))g, \pi ( u+\frac{1}{2}j\inv   (v))g \Big\rangle \Big| \\
&\leq & H(j\inv (v)) \, \quad \quad \forall u\in \rdd \, .
\end{eqnarray*}
We conclude that
\begin{equation}
  \label{eq:b1}
\intrdd \sup _{u\in \rdd }   |\cV _\Phi \sigma (u,\zeta)| \, \,
v(j\inv (\zeta )) \, d\zeta  \leq \intrdd H(v\inv (\zeta )) v(j\inv
(\zeta )) \, d\zeta = \|H\|_{L^1_v}  < \infty  \, ,
  \end{equation}
and so  $\sigma \in \mif _{\vv } (\rdd )$. 

\vspace{3 mm}

The discrete condition (iii) is similar, but technically more subtle
to handle. 

(i) $\, \Longrightarrow \, $ (iii) $\quad $ To show this
implication, we use the well-known fact that the \stft\ of a
distribution possesses ``nice'' local properties, see~\cite{CG03,book,Toft04} for
various  statements and proofs. In particular, if $\sigma \in \mif
_{\vv }$, then $\cV _\Phi \sigma \in W(C, \ell ^{\infty ,1}_{\vv }
)(\fourd )$~\cite[Thm.~12.2.1]{book}. This means the following: let
$Q=[-1/2,1/2]^{2d}$ and define the sequence $a_k, k\in \zdd, $ to be 
$a_k = \sup _{\zeta\in k+Q} \sup _{z\in \rdd } |\cV _\Phi \sigma (z,\zeta)| $;
then $$\sum _{k\in \zdd } a_k \,  v(j\inv (k)) = \|a\|_{\ell ^1 _{\vv}}
\leq C \|\sigma \|_{\mif _{\vv}} < \infty \, .$$
Using~\eqref{eq:4a} once more, we obtain that
$$
|\langle \sigma ^w \pi (\mu )g, \pi (\lambda )g\rangle | = \Big|\cV _\Phi
\sigma \big( \frac{\lambda +\mu }{2}, j(\lambda -\mu)\big) \Big| \leq a_k
\, \quad \quad \mathrm{ if } \,\,  j(\lambda -\mu ) \in k+Q \, .
$$
 Now set
\begin{equation*}
  \label{eq:11}
  h(\lambda ) = a_k \quad \quad \text{ if } \lambda \in j\inv (k+Q) =
  j\inv (k)+Q \, .
\end{equation*}
Then
\begin{eqnarray}
  \sum _{\lambda \in \Lambda } h(\lambda ) \, v(\lambda ) &=& \sum _{k\in
    \zdd } \sum _{\lambda \in j\inv (k) +Q} a_k \,\,   v(\lambda )
  \notag \\
&\leq & \sum _{k\in     \zdd } \sum _{\lambda \in j\inv (k) +Q} a_k \,
v(j\inv (k )) \, \sup _{u\in Q} v(-u) \label{eq:b2} \\
&=& C \max _{k\in \zdd } \mathrm{card}\, \{\lambda \in \Lambda :
\lambda \in j\inv (k) +Q\}| \, 
\sum _{k\in \zdd } a_k \,  v(j\inv (k)) \notag  \\
& \leq & C' \|\sigma \|_{\mif _{\vv }}
\, . \notag 
\end{eqnarray}
This is (iii) as desired. 

(iii) $\, \Longrightarrow \, $ (ii) $\quad $ To prove this implication,
we finally use the hypothesis that $\cG (g,\Lambda )$ is a Gabor frame. Since
$g\in M^1_{v }$, the dual window $\gamma $ is also in $M^1_{v }$
by~Theorem~\ref{tfmain}. In particular, every \tfs\ $\pi (u)g$ has the
following frame expansion:
\begin{equation}
  \label{eq:12}
\pi (u)g = \sum _{\nu \in \Lambda } \langle \pi (u)g,\pi (\nu
) \gamma \rangle \pi (\nu ) g \, .
\end{equation}
If $g,\gamma \in \mvv $, then by the local properties of \stft
s~\cite[Thm.~12.2.1]{book}, we 
know that $V_\gamma g \in 
W(C,\ell ^1_{v })(\rdd )$. This  means that for every
relatively compact set $C\subseteq \rdd $ we have
$$
\sum _{\nu \in \Lambda } \sup _{u\in C} |V_\gamma g(\nu +u)|
v(\nu ) \leq C \|g\|_{M^1_{v }}
$$
In particular, if  $C$ is  a relatively compact fundamental domain of
the lattice $\Lambda $  and  
\begin{equation}
  \label{eq:13b}
\alpha (\nu ) = \sup _{u\in C} |V_\gamma g(\nu +u)| = \sup _{u\in C}
|\langle \pi (-u)g,\pi (\nu) \gamma \rangle | \, ,
\end{equation}
then the sequence $\alpha $ is in $\ell ^1 _{v }(\Lambda )$.

Given $z,w \in \rdd $ we can write them uniquely as $w=\lambda + u,
z=\mu +u'$ for  $\lambda ,\mu \in \Lambda$ and $u,u' \in
C$. Inserting  the expansions~\eqref{eq:12} and the definition of $\alpha$
in the matrix entries,  we find that
\begin{eqnarray*}
\lefteqn{  |\langle \sigma ^w \pi (\mu +u')g, \pi (\lambda +u)g\rangle | =
  |\langle \sigma ^w \pi (\mu ) \pi (u')g, \pi (\lambda) \pi
  (u)g\rangle |} \\ 
&\leq & \sum _{\nu, \nu' \in \Lambda } |\langle \sigma ^w \pi (\mu
+\nu ')g, \pi (\lambda +\nu )g\rangle| \, |\langle \pi (u')g, \pi
(\nu') \gamma \rangle | \, |\langle \pi (u) g, \pi (\nu )\gamma
\rangle | \\
&\leq & \sum _{\nu, \nu '\in \Lambda } h(\lambda +\nu - \mu -\nu ')
\alpha (\nu') \alpha (\nu ) \\
&=& (h\ast \alpha \ast \check{\alpha }  )(\lambda -\mu )\, ,
\end{eqnarray*}
with $\check{\alpha }(\lambda ) = \alpha (-\lambda )$. 
Since  $h \in  \ell ^1 _{v }$ by hypothesis (iii)  and $\alpha
\in \ell ^1 _{v }$ by construction, we
also have $h\ast \alpha \ast \check{\alpha} \in \ell ^1 _{v
}(\Lambda )$. 

Now set 
$$
H(z) = \sum _{\lambda \in \Lambda } (h\ast \alpha  \ast \check{\alpha }
) (\lambda ) \,  \chi _{C-C}(z-\lambda) \, .
$$
Then 
\begin{equation}
  \label{eq:b3}
 \|H\|_{L^1_{v}} \leq \sum _{\lambda }  h\ast \alpha  \ast
\check{\alpha }
 (\lambda ) v(\lambda ) \, \, \|\chi _{C-C}\|_{L^1_v} = c \|h\ast
\alpha \ast \check{\alpha }  \|_{\ell ^1_v }<\infty \, . 
\end{equation}
If 
$z,w \in \rdd $ with  $w=\lambda + u, z=\mu +u'$ for $\lambda ,\mu \in
\Lambda$ and $u,u'\in C$, then $w-z\in \lambda -\mu +C-C$ and 
$(h\ast \alpha  \ast \alpha
^*) (\lambda -\mu ) \leq  H(w-z)$. Combining these observations, we
have shown that
$$\Big| \langle \sigma ^w \pi (z)g, \pi (w)g\rangle \Big| \leq (h\ast \alpha \ast
\check{\alpha })(\lambda -\mu ) \leq  H(w-z)\, ,
$$
and this is (ii). 
\end{proof}

\begin{cor}
  Under the hypotheses of Theorem~\ref{almost}, assume that $T:\cS (\rd )
  \to \cS ' (\rd )$ is continuous and    satisfies the estimates
$$
\Big| \langle T \pi (\mu )g, \pi (\lambda )g\rangle \Big| \leq
h(\lambda -\mu ) \quad \quad \forall \lambda,\mu \in \Lambda 
$$
for some $h\in \ell ^1_v$. Then $T=\sigma ^w $ for some   symbol $\sigma \in \mif
_{\vv}$.
\end{cor}

\begin{proof}
  Schwartz's kernel theorem and \eqref{eq:f2} imply   that $T=
 \sigma ^w $ for some distributional symbol $\sigma \in \cS ' (\rdd
 )$ (see also~\cite[Thm.~14.3.5]{book}). Now apply Theorem~\ref{almost}. 
\end{proof}
\rems\ 1. 
 Motivated by the concept of ``confined symbols''~\cite{BL89}, \sjo\ proved that
 $\sigma \in \mif $ \fif\ there 
exists $h\in \ell ^1(\Lambda )$ such that \\ $\|(T_\mu \chi  )^w \weyl
(T_\lambda \chi )^w \|_{L^2 \to L^2} \leq h(\lambda -\mu )$, where
$\chi \in \cS (\rdd )$ satisfies $\sum_{\lambda \in \Lambda } \chi (t-\lambda
) = 1 $. The equivalence (i) $\, \Longleftrightarrow \, $ (ii)  was
also obtained independently by  Strohmer\cite{strohmer04}.

2. Property (ii) says that $\sigma ^w $ preserves the \tf\
localization and that $\sigma ^w$ maps  the \tfs s $\pi (z)g$
into functions in $\mvv $  with a uniform  envelope $H$  in the \tf\ plane. 
This could be rephrased by saying that  $\sigma ^w$ maps \tf\ ``atoms'' into
\tf\ ``molecules''. 

3. By property (iii)  $\sigma ^w$ is almost diagonalized by the
Gabor frame $\cG (g,\Lambda )$. 
It is well-known that certain types of \psdo s are almost diagonalized
with respect to wavelet bases or local Fourier
bases~\cite{meyer-2,tachizawa-rochberg98}. What is
remarkable in Theorem~\ref{almost}  is that the almost diagonalization
property actually characterizes a symbol class.

\subsection{Matrix Formulation}

Let us formulate Theorem~\ref{almost} on a more conceptual level. 
Let $f = \sum _{\mu \in \Lambda } \langle f, \pi (\mu
  )\gamma\rangle \pi (\mu ) g$ be the Gabor expansion of $f\in \lrd
$, then 
\begin{equation}
  \label{eq:18}
  C_g (\sigma ^wf)(\lambda ) = \langle \weyl f, \pi (\lambda )g\rangle
  = \sum _{\mu \in \Lambda } \langle f, \pi (\mu 
  )\gamma\rangle \, \langle \weyl \pi (\mu )g,  \pi (\lambda ) g
  \rangle \, .
\end{equation}
We therefore define the  matrix $M(\sigma )$ associated  to the 
symbol  $\sigma $  with respect to a  Gabor frame by the   entries 
\begin{equation}
  \label{eq:11b}
  M(\sigma ) _{\lambda \mu } = \langle \sigma ^w \pi (\mu )g, \pi
  (\lambda )g\rangle \, , \quad \quad \lambda, \mu \in \Lambda \, .
\end{equation}
With this notation,  \eqref{eq:18} can be recast as 
\begin{equation}
  \label{eq:19}
  C_g(\weyl f) = M(\sigma) C_\gamma f \, ;
\end{equation}
or as a commutative  diagram:
\begin{equation} \label{diagram1}
\begin{matrix}
&\lrd   & \stackrel{\sigma ^w}{\longrightarrow} & \lrd  & \cr
&\downarrow C_\gamma & &\downarrow C_g& \cr &\ell ^2(\Lambda)\, 
&\stackrel{M(\sigma )}{\longrightarrow} &\ell ^2(\Lambda ) &
\end{matrix}
\end{equation}


\begin{lemma} \label{matrixprop}
  If $\sigma ^w $ is bounded on $\lrd $, then $M(\sigma )
  $ is bounded on $\ell ^2( \Lambda )$ and maps $\mathrm{ran}\, C_g$
  into $\mathrm{ran}\, C_g$ with $\mathrm{ker}\, M(\sigma ) \supseteq
  (\mathrm{ran}\, C_g )^\perp= \mathrm{ker} \, C_g ^*$. 
\end{lemma}

\begin{proof}
Note that $\mathrm{ran} \, C_\gamma = \mathrm{ran} \, C_g$, since
$\langle f, \pi (\lambda ) \gamma \rangle = \langle f, \pi (\lambda )
S\inv g \rangle = \langle S\inv f, \pi (\lambda ) g \rangle $ for all
$\lambda \in \Lambda $,  or $C_\gamma = C_g S\inv $.   

Consequently,   by  the frame property and ~\eqref{diagram1} we have 
$$
\|M(\sigma ) C_\gamma f \|_2 = \|C_g (\weyl f )\|_2 \leq C_1 \|\weyl f
\|_2 \leq C_2 \|f\|_2 \leq C_3 \|C_g f\|_2\, , 
$$
and so $M(\sigma ) $ is bounded from  $\mathrm{ran}\, C_g$
  into $\mathrm{ran}\, C_g$. If $c\in (\mathrm{ran}\, C_g )^\perp=
  \mathrm{ker} \, C_g ^*$, then $\sum _{\mu \in \Lambda } c_\mu \pi
  (\mu )g = 0$, and thus $(M(\sigma ) c)(\lambda ) = \sum _{\mu \in
    \Lambda } \langle \weyl \pi (\mu ) g, \pi (\lambda ) g\rangle
  c_\mu = 0$, i.e., $c\in  \mathrm{ker}\, M(\sigma )$. 
\end{proof}

Since  $\cG (g, \Lambda )= \{\pi (\lambda )g: \lambda \in \Lambda
\}$ is only a frame, but not a basis,  not every matrix $A$ is of the
form $M(\sigma )$. 
It is easy to see that the  properties of Lemma~\ref{matrixprop} imply
that  $A = M(\sigma )$ for some $\sigma \in \cS ' (\rdd )$. 


Next we formalize the properties of the matrices occurring in
Theorem~\ref{almost}.

\begin{definition}
We say that a matrix $A=(a_{\lambda \mu})_{\lambda, \mu \in \Lambda}$
belongs to  $\cC _v = \cC _v (\Lambda )$, if there exists a sequence
$h \in \ell ^1 _v (\Lambda )$ such that 
\begin{equation}
  \label{eq:8}
  |a_{\lambda \mu }| \leq h(\lambda -\mu ) \quad \quad \forall \lambda
  ,  \mu \in \Lambda \, .
\end{equation}
We endow  $\cC _v$  with the norm 
\begin{eqnarray}
  \|A\|_{\cC _v} &=& \inf \{ \|h\|_{\ell ^1_v}:   |a_{\lambda \mu }|
  \leq h(\lambda -\mu ),  \forall \lambda    ,  \mu \in \Lambda \}
  \label{eq:10} \\
&=& \sum _{\mu \in \Lambda } \sup _{\lambda \in \Lambda }  |a_{\lambda
  , \lambda -\mu }|\,  v(\mu ) \, .   
\notag 
  \end{eqnarray}
\end{definition}

Since every $A\in \cC _v$ is dominated by a convolution operator, the
algebra property is evident. 

\begin{lemma}
  $\cC _v $ is a Banach $*$-algebra. 
\end{lemma}

\rem\   If $A\in \cC _v$, then $A$ is automatically bounded on
$\ell ^p_m $ for $1\leq p \leq \infty $ and $m \in \mv $. This follows
from the pointwise inequality $|Ac(\lambda )| \leq (h\ast |c|)(\lambda
)$ and Young's inequality. If $\Lambda = \alpha \zd \times \beta \zd $, then also
\begin{equation}
  \label{eq:20}
  \|Ac\|_{\lpqm } \leq \|h\ast |c|\, \|_{\lpqm } \leq \|h\|_{\ell ^1_v
    } \,  \|c\|_{\lpqm } \, .
\end{equation}

Theorem~\ref{almost} can be recast as follows.

\begin{tm} \label{almost2}
  A symbol $\sigma $ is in $\mif _{\vv } $ \fif\ $M(\sigma ) \in \cC
  _v$ and 
  \begin{equation}
    \label{eq:23}
\| \sigma \|_{\mif _{\vv }} \asymp     \|M(\sigma ) \|_{\cC _v} \, .
  \end{equation}
\end{tm}
The estimate $\|M(\sigma ) \|_{\cC _v} \leq C_1 \|\sigma \|_{\mif
  _{\vv }}$ is contained in \eqref{aha}, the converse inequality follows
  by combining~\eqref{eq:b1} and \eqref{eq:b3}. 

\vspace{4 mm}

\rem\  In view of this  reformulation  it is natural to consider other
matrix algebras and study the relation between symbols and the
membership of  $M(\sigma )$  in a  matrix algebra.

\section{The Proofs of \sjo 's Results}

We are now ready to prove \sjo 's results in their ``natural'' context
and at the same time we formulate suitable extensions. In place of
``hard analysis'' we use  \tf\ methods,  Theorem~\ref{almost},  and
recent Banach algebra techniques. 

Though frames do not enter in the formulation of the results, they are
vital in the proofs. To treat all weights in the class $\mv $, we need
to assume as in Theorem~\ref{tfmain} that the window is chosen from an
appropriate space of test functions $\mvv $.

\subsection{Boundedness}

\begin{tm} \label{bound}
If $\sigma \in \mif _{\vv }$, then $\sigma ^w $ is bounded on $\Mmpq $
for $1\leq p,q \leq \infty  $ and all $m\in \cM _{v }$. The operator
norm can be estimated uniformly by 
$$
\|\weyl \|_{\Mmpq \to \Mmpq} \leq C \|M(\sigma )\|_{\cC _v} \asymp
\|\sigma \|_{\mif _{\vv }} \, ,
$$
with a  constant  independent of $p,q,$ and $m$. 
\end{tm}

\begin{proof}
  Fix a Gabor frame $\cG (g, \alpha \zd \times \beta \zd )$ with
  window $g \in M^1_{\vv }$. By
  Theorem~\ref{tfmain}  also $\gamma \in M^1_{\vv }$ and  the
  following  norms are equivalent  on
  $\Mmpq $: $   \|f\|_{\Mmpq} \asymp \|C_g f\|_{\lpqm  }\asymp \|C_\gamma f\|_{\lpqm
  }$ for every $1\leq p,q \leq \infty$ and $m\in \cM _{v }$. 

Now let $f \in M^1_{v }\subseteq \lrd $ be arbitrary. Applying
diagram~\eqref{diagram1}, we  estimate the $\Mmpq $-norm 
of $\weyl f$ as follows:
$$
 \| \sigma ^w f \|_{M^{p,q}_m} \leq C_0 \,  \| C_\gamma (\sigma ^w f)
 \|_{\ell    ^{p,q}_m} =  C_0\,   \| M(\sigma ) C_gf \| _{\ell
   ^{p,q}_m} \, .  
$$
Since $M(\sigma ) \in \cC _v$ by Theorem~\ref{almost}, $M(\sigma )$ is
bounded on $\lpqm $ for $m\in \cM _v$ by~\eqref{eq:20}. So  we
continue the above estimate by  
$$
\| \sigma ^w f \|_{M^{p,q}_m} \leq C_0 \|M(\sigma )\|_{\lpqm \to \lpqm
} \, \, \|C_g f\|_{\lpqm } \leq C_1 \|M(\sigma )\|_{\cC _v}\,  \|f\|_{\Mmpq } \, .
$$
This implies that $\weyl $ is bounded on the closure of $M^1_{v }$
in the $\Mmpq $-norm. If $p,q<\infty$, then by density 
$\weyl $ is bounded on $\Mmpq$. For $p=\infty $ or $q=\infty $, the
argument has to be modified as in~\cite{BGHO04}.
\end{proof}

\rems\ 1. In particular, if $\sigma \in \mif $, then $\weyl $ is
bounded on $\lrd $~\cite{Sjo94,Bou97} and on all  $M^{p,q}(\rd )$ for
$1\leq p,q \leq \infty$~\cite{GH99,book}. 

2.  Theorem~\ref{bound} is a slight improvement
over~\cite[Thm.~14.5.6]{book} where the boundedness on $\Mmpq $ for
$m\in \cM _{v } $ required that $\sigma \in \mif _w$ with
$w(\zeta) = v(j\inv(\zeta ) /2)^2 \geq v(j\inv(\zeta ) )$. 

Since $S^0_{0,0} \subseteq \mif $, the Weyl transforms $\weyl $ for
$\sigma \in \mif $ cannot be bounded on $L^p (\rd )$ in general. Using
the embeddings $ L^p \subseteq M^{p,p'}$ for $1\leq
p\leq 2$ and $ L^p \subseteq M^p$ for $2\leq p \leq
\infty $, we obtain an $L^p$ result as follows.

\begin{cor}
  Assume that $\sigma \in \mif $. If $1\leq p\leq 2$, then $\weyl $
  maps $L^p $ into $M^{p,p'}$, whereas for  $2\leq p\leq \infty $,  $\weyl$
  maps $L^p $ into $M^p$. 
\end{cor}

\subsection{The Algebra Property}

\begin{tm}\label{sjowei}
If $v$ is submultiplicative, then $\mif _v $ is a Banach $*$-algebra
with respect to the 
twisted product $\sharp$ and the involution $\sigma \to \bar{\sigma}$.
\end{tm}

\begin{proof} 
It is convenient to use a tight Gabor frame $\cG (g, \alpha \zd \times
\beta \zd )$  with   $\gamma =  g\in \mvv $ as in \eqref{eq:16b}. By
using~\eqref{diagram1} twice, we obtain that
\begin{eqnarray*}
  M(\sigma \, \sharp \, \tau ) C_g f &=& C_g ( (\sigma \, \sharp \,
  \tau )^w  \, f) = C_g (\weyl \, \tau ^w f) \\
&=& M(\sigma ) \big( C_g (\tau ^w f )) = M(\sigma ) M(\tau ) C_g f \, . 
\end{eqnarray*}
Therefore the operators $M(\sigma \, \sharp \tau )$ and $M(\sigma )
M(\tau )$ coincide on $\mathrm{ran} \, C_g$. 
Since  \\
$M(\sigma ) | _{(\mathrm{ran} \, C_g )^\perp } = 0$ for
all $\sigma \in \mif $ by Lemma~\ref{matrixprop}, we obtain that 
\begin{equation}
  \label{eq:21}
M(\sigma \, \sharp \, \tau ) =   M(\sigma ) M(\tau )
\end{equation}
as an identity of matrices (on $\ell ^2$). 

Now, if  $\sigma , \tau \in \mif _v$, then $M(\sigma ), M(\tau ) \in
\cC _{v\circ j}$ by Theorem ~\ref{almost2}. By the algebra property of
$\cC _{v\circ j}$  we have $M(\sigma ) M(\tau ) \in \cC _{v\circ j}$,
and once again by Theorem~\ref{almost2} we have $M(\sigma \, \sharp \,
\tau ) \in \cC _{v\circ j}$ with the norm estimate
$$
\|\sigma \, \sharp \, \tau  \| _{\mif _v} \leq C_0 \|M(\sigma \, \sharp \,
\tau )\|_{\cC _{v\circ j}} \leq  C_0 \, \|M(\sigma )\|_{\cC _{v\circ j}}\,
\|M(\tau )\|_{\cC _{v\circ j}} \leq C_1 \, \|\sigma \|_{\mif _v} \,
\|\tau \|_{\mif _v} \, . 
$$
\end{proof}
Compare~\cite{grocomp,Sjo94,Sjo95,Toft01} for other proofs. 

\subsection{Wiener Property of Sj\"ostrand's Class}

For the Wiener property we start with two results about
the Banach algebra $\cC _v$.

\begin{tm} \label{bask}
Assume that  $v$ is a submultiplicative weight  satisfying the
GRS-condition 
\begin{equation}
  \label{eq:ne8}
\lim _{n\to \infty } v(nz)^{1/n} = 1 \quad \forall z\in \rdd \, .  
\end{equation}
  If $A\in \cC _v$ and $A$ is invertible on $\ell ^2(\zd )$, then 
$A\inv \in \cC _v$. As a consequence
\begin{equation}
  \label{eq:14b}
\mathrm{Sp}_{\cB (\ell ^2)} (A) = \mathrm{Sp}_{\cC _v } (A)   
\end{equation}
for all $A \in \cC _v$, where $\mathrm{Sp}_{\cA } (A)$ denotes the
spectrum of $A$ in the algebra $\cA $.  
\end{tm}

Originally, this important result  was proved by
Baskakov~\cite{Bas90,Bas97a}  in several papers, and by
Sj\"ostrand~\cite{Sjo95}   for the   
unweighted  case $v \equiv 1$. 

Recall that an operator $A: \ell ^2 \to \ell ^2$ is pseudo-invertible,
if there exists a closed subspace $\cR \subseteq \ell ^2$, such
that $A$ is invertible on $\cR$ and $\mathrm{ker}\, A = \cR
^\perp $. The unique operator $A^\dagger $ that satisfies $A^\dagger
A h = A A^\dagger h = h $ for $h\in \cR$ and $\mathrm{ker}\,
A^\dagger = \cR^\perp $ is called the (Moore-Penrose)
pseudo-inverse of $A$. The following lemma
is borrowed from ~\cite{FoG04}.

  \begin{lemma}[Pseudoinverses]   \label{pseudo}
    If $A\in \cC _v$ has a (Moore-Penrose) pseudoinverse $A^\dagger$,
then $A^\dagger \in \cC _v $. 
  \end{lemma}

  \begin{proof}
By means of the Riesz functional calculus~\cite{rudin73}  the
pseudoinverse can be written as 
$$
A^{\dagger }  = \frac{1}{2\pi i} \int _{ C } \frac{1}{z}\,
  (z\mathrm{I}-A)\inv \, dz \, ,
$$
where $C$ is a suitable path surrounding
$\mathrm{Sp}_{\cB (\ell ^2)} (A) \setminus \{ 0\}$.   By
\eqref{eq:14b} this formula   
make sense in  $\cC _v$,  and consequently  $A^\dagger \in \cC _v $.     
  \end{proof}

\begin{tm}\label{wiener}
 Assume that $v$ satisfies the GRS-condition $\lim _{n\to \infty }
 v(nx)^{1/n}  =1, \forall   x\in \rdd $. 
If   $\sigma \in M^{\infty , 1}_v (\rdd )$ and  $\sigma
  ^w $ is  invertible on $\lrd $, 
then $(\sigma ^w) \inv = \tau ^w  $ for some $\tau \in M^{\infty ,1}_v$.  
\end{tm}

\begin{proof}
Again, we use a tight Gabor frame $\cG (g, \alpha \zd \times \beta \zd
)$ with $g = \gamma \in \mvv $ as in~\eqref{eq:16b} for the analysis
of the Weyl transform. 
 
Let $\tau \in \cS ' (\rdd )$ be the unique distribution 
such that $\tau ^w = (\weyl )\inv $. Then the  matrix $M(\tau ) $ is
 bounded on $\ell ^2 $ and  maps
$\mathrm{ran} \, C_g $ into $\mathrm{ran} \, C_g $ with  $\mathrm{ker}
\, M(\tau ) \subseteq (\mathrm{ran}\, C_g ) ^{\perp }$ (by
Lemma~\ref{matrixprop}). 

We show that $M(\tau )$ is the pseudo-inverse of $M(\sigma )$. 
Let $c= C_gf \in \mathrm{ran}\, C_g$, then 
$$
  M(\tau ) M(\sigma ) C_g f = M(\tau ) C_g (\weyl f) =   C_g (\tau
  ^w \weyl f) = C_g f \, , 
$$
and likewise $M(\sigma ) M(\tau ) = \mathrm{I}_{\mathrm{ran}\,
  C_g}$. Since  $\mathrm{ker}\, M(\sigma ),\mathrm{ker}\, M(\tau )
\subseteq   (\mathrm{ran} \, C_g )^{\perp} 
$,  we  conclude that $M(\tau ) = M(\sigma ) ^\dagger $. 

By Theorem~\ref{almost} the hypothesis  $\sigma \in \mif _v $ implies
that $M(\sigma )$ belongs to the matrix algebra $\cC _{v\circ
  j}$. Consequently by Lemma~\ref{pseudo}, we also have $M(\tau ) = M(\sigma )^\dagger \in
\cC _{v\circ j}$. Using Theorem~\ref{almost} again,  we conclude that
$\tau \in \mif _v$. This finishes the proof of the Wiener property. 
\end{proof}

It can be shown that Theorem~\ref{wiener} is false, when $v$ does not
satisfy~\eqref{eq:ne8}. Thus the  GRS-condition is sharp. 

\vspace{1 cm}

\begin{cor}
[Spectral Invariance on Modulation Spaces] 
If  $\sigma \in \mif _{\vv }$ and   $\sigma ^w $ is invertible  on $\lrd $,
then   $\sigma ^w$ 
is invertible simulaneously  on all \modsp s $M^{p,q}_m(\rd )$, where
$1\leq p,q\leq \infty $ and  $m\in \cM _{v }$. 
\end{cor}

\begin{proof}
By Theorem~\ref{wiener} $(\weyl )\inv = \tau ^w$ for some $\tau \in \mif
_{\vv }$ and then  by Theorem~\ref{bound} $\tau ^w$ is bounded on $\Mmpq $
for the range of $p,q$ and $m$ specified. Since $\weyl \tau ^w = \tau
^w \weyl = \mathrm{I}$ on $\mvv $, this factorization extends by
density to all of $\Mmpq$. Thus $\tau ^w = (\weyl )\inv $ on $\Mmpq
$. 
\end{proof}

\rems\ 1. It is known that $\mif $ is invariant under convolution with
``chirps'' $e^{it\cdot Ct}$ for any symmetric real-valued $d\times
d$-matrix $C$~\cite{Sjo94,book}. As a consequence,  the properties of
the symbol class $\mif $ carry over to other calculi of \psdo s, in
particular to the \knc ~\cite{Sjo94,Toft01}.

2. Translation and modulation operators can be defined on arbitrary
locally compact abelian groups (LCA groups), and consequently,  \modsp s and
the \knc\ are 
well-defined on LCA groups  in place of $\rd $. 
Therefore \sjo 's results should hold in the general context of LCA
groups, but it is clear that the methods of classical analysis can no
longer be  applied, whereas it is plausible that \tf\ methods can be
generalized. For instance, it is not hard to verify that the matrix
algebra $\cC _v $ for $v \equiv 1$ coincides with $\mif (\zd \times
\bT ^d)$. Thus Theorem~\ref{bask} says that the Wiener property holds for
the \modsp\ $\cC = \mif (\zd \times \bT ^d)$. Therefore we conjecture that
Theorem~\ref{wiener} holds for $\mif (\cG \times \widehat{\cG })$ for
an arbitrary LCA group $\cG $, and will pursue this question
elsewhere.




\def\cprime{$'$}



\begin{thebibliography}{10}

\bibitem{Bas90}
A.~G. Baskakov.
\newblock Wiener's theorem and asymptotic estimates for elements of inverse
  matrices.
\newblock {\em Funktsional. Anal. i Prilozhen.}, 24(3):64--65, 1990.

\bibitem{Bas97a}
A.~G. Baskakov.
\newblock Asymptotic estimates for elements of matrices of inverse operators,
  and harmonic analysis.
\newblock {\em Sibirsk. Mat. Zh.}, 38(1):14--28, i, 1997.

\bibitem{beals77}
R.~Beals.
\newblock Characterization of pseudodifferential operators and applications.
\newblock {\em Duke Math. J.}, 44(1):45--57, 1977.

\bibitem{BHW95}
J.~J. Benedetto, C.~Heil, and D.~F. Walnut.
\newblock Differentiation and the {B}alian--{L}ow theorem.
\newblock {\em J. Fourier Anal. Appl.}, 1(4):355--402, 1995.

\bibitem{BGHO04}
A.~Benyi, K.~Gr{\"o}chenig, C.~Heil, and K.~Okoudjou.
\newblock Modulation spaces and a class of bounded multilinear
  pseudodifferential operators.
\newblock {\em J. Operator Theory}, 2004.
\newblock To appear.

\bibitem{bjork66}
G.~Bj{\"o}rk.
\newblock Linear partial differential operators and generalized distributions.
\newblock {\em Ark. Mat.}, 6:351--407, 1966.

\bibitem{BCG04}
P.~Boggiatto, E.~Cordero, and K.~Gr{\"o}chenig.
\newblock Generalized anti-{W}ick operators with symbols in distributional
  {S}obolev spaces.
\newblock {\em Integral Equations Operator Theory}, 48(4):427--442, 2004.

\bibitem{BL89}
J.-M. Bony and N.~Lerner.
\newblock Quantification asymptotique et microlocalisations d'ordre
  sup\'erieur. {I}.
\newblock {\em Ann. Sci. \'Ecole Norm. Sup. (4)}, 22(3):377--433, 1989.

\bibitem{Bou97}
A.~Boulkhemair.
\newblock Remarks on a {W}iener type pseudodifferential algebra and {F}ourier
  integral operators.
\newblock {\em Math. Res. Lett.}, 4(1):53--67, 1997.

\bibitem{Boul99}
A.~Boulkhemair.
\newblock {$L\sp 2$} estimates for {W}eyl quantization.
\newblock {\em J. Funct. Anal.}, 165(1):173--204, 1999.

\bibitem{CG03}
E.~Cordero and K.~Gr{\"o}chenig.
\newblock Time-frequency analysis of localization operators.
\newblock {\em J. Funct. Anal.}, 205(1):107--131, 2003.

\bibitem{CR02}
W.~Czaja and Z.~Rzeszotnik.
\newblock Pseudodifferential operators and {G}abor frames: spectral
  asymptotics.
\newblock {\em Math. Nachr.}, 233/234:77--88, 2002.

\bibitem{daubechies90}
I.~Daubechies.
\newblock The wavelet transform, time-frequency localization and signal
  analysis.
\newblock {\em IEEE Trans. Inform. Theory}, 36(5):961--1005, 1990.

\bibitem{feichtinger81}
H.~G. Feichtinger.
\newblock On a new {S}egal algebra.
\newblock {\em Monatsh. Math.}, 92(4):269--289, 1981.

\bibitem{feichtinger-modulation1}
H.~G. Feichtinger.
\newblock Modulation spaces on locally compact abelian groups.
\newblock Technical report, University of Vienna, 1983.

\bibitem{fei83}
H.~G. Feichtinger.
\newblock Modulation spaces on locally compact abelian groups.
\newblock In {\em Proceedings of ``International Conference on Wavelets and
  Applications" 2002}, pages 99--140, Chennai, India, 2003.
\newblock Updated version of a technical report, University of Vienna, 1983.

\bibitem{fg97jfa}
H.~G. Feichtinger and K.~Gr{\"o}chenig.
\newblock Gabor frames and time-frequency analysis of distributions.
\newblock {\em J. Functional Anal.}, 146(2):464--495, 1997.

\bibitem{folland89}
G.~B. Folland.
\newblock {\em Harmonic Analysis in Phase Space}.
\newblock Princeton Univ. Press, Princeton, NJ, 1989.

\bibitem{FoG04}
M.~Fornasier and K.~Gr{\"o}chenig.
\newblock Intrinsic localization of frames.
\newblock {\em Preprint}, 2004.

\bibitem{gelfandraikov}
I.~Gel'fand, D.~Raikov, and G.~Shilov.
\newblock {\em Commutative normed rings}.
\newblock Chelsea Publishing Co., New York, 1964.

\bibitem{book}
K.~Gr{\"o}chenig.
\newblock {\em Foundations of time-frequency analysis}.
\newblock Birkh\"auser Boston Inc., Boston, MA, 2001.

\bibitem{grocomp}
K.~Gr{\"o}chenig.
\newblock Composition and spectral invariance of pseudodifferential operators
  on modulation spaces.
\newblock {\em Preprint}, 2004.

\bibitem{GH99}
K.~Gr{\"o}chenig and C.~Heil.
\newblock Modulation spaces and pseudodifferential operators.
\newblock {\em Integral Equations Operator Theory}, 34(4):439--457, 1999.

\bibitem{GH03}
K.~Gr{\"o}chenig and C.~Heil.
\newblock Modulation spaces as symbol classes for pseudodifferential operators.
\newblock In S.~T. M.~Krishna, R.~Radha, editor, {\em Wavelets and Their
  Applications}, pages 151--170. Allied Publishers, Chennai, 2003.

\bibitem{GH04}
K.~Gr{\"o}chenig and C.~Heil.
\newblock Counterexamples for boundedness of pseudodifferential operators.
\newblock {\em Osaka J. Math}, 41:1--11, 2004.

\bibitem{GL03}
K.~Gr{\"o}chenig and M.~Leinert.
\newblock Wiener's lemma for twisted convolution and {G}abor frames.
\newblock {\em J. Amer. Math. Soc.}, 17:1--18, 2004.

\bibitem{heil-ramanathan97}
C.~Heil, J.~Ramanathan, and P.~Topiwala.
\newblock Singular values of compact pseudodifferential operators.
\newblock {\em J. Functional Anal.}, 150(2):426--452, 1997.

\bibitem{her01}
F.~H{\'e}rau.
\newblock Melin-{H}\"ormander inequality in a {W}iener type pseudo-differential
  algebra.
\newblock {\em Ark. Mat.}, 39(2):311--338, 2001.

\bibitem{hormander3}
L.~H{\"o}rmander.
\newblock {\em The analysis of linear partial differential operators. {III}},
  volume 274 of {\em Grundlehren der Mathematischen Wissenschaften [Fundamental
  Principles of Mathematical Sciences]}.
\newblock Springer-Verlag, Berlin, 1994.
\newblock Pseudo-differential operators, Corrected reprint of the 1985
  original.

\bibitem{janssen95}
A.~J. E.~M. Janssen.
\newblock Duality and biorthogonality for {W}eyl-{H}eisenberg frames.
\newblock {\em J. Fourier Anal. Appl.}, 1(4):403--436, 1995.

\bibitem{lab01}
D.~Labate.
\newblock Pseudodifferential operators on modulation spaces.
\newblock {\em J. Math. Anal. Appl.}, 262(1):242--255, 2001.

\bibitem{lab02}
D.~Labate.
\newblock Time-frequency analysis of pseudodifferential operators.
\newblock {\em Monatsh. Math.}, 133(2):143--156, 2001.

\bibitem{meyer-2}
Y.~Meyer.
\newblock {\em Ondelettes et op\'erateurs. {I}{I}}.
\newblock Hermann, Paris, 1990.
\newblock Op\'erateurs de Calder\'on-Zygmund. [Calder\'on-Zygmund operators].

\bibitem{PT04}
S.~Pilipovi{\'c} and N.~Teofanov.
\newblock Pseudodifferential operators on ultra-modulation spaces.
\newblock {\em J. Funct. Anal.}, 208(1):194--228, 2004.

\bibitem{tachizawa-rochberg98}
R.~Rochberg and K.~Tachizawa.
\newblock Pseudodifferential operators, {G}abor frames, and local trigonometric
  bases.
\newblock In {\em Gabor analysis and algorithms}, pages 171--192. Birkh\"auser
  Boston, Boston, MA, 1998.

\bibitem{rudin73}
W.~Rudin.
\newblock {\em Functional analysis}.
\newblock McGraw-Hill Book Co., New York, 1973.
\newblock McGraw-Hill Series in Higher Mathematics.

\bibitem{shubin91}
M.~A. Shubin.
\newblock {\em Pseudodifferential Operators and Spectral Theory}.
\newblock Springer-Verlag, Berlin, second edition, 2001.
\newblock Translated from the 1978 Russian original by Stig I. Andersson.

\bibitem{Sjo94}
J.~Sj{\"o}strand.
\newblock An algebra of pseudodifferential operators.
\newblock {\em Math. Res. Lett.}, 1(2):185--192, 1994.

\bibitem{Sjo95}
J.~Sj{\"o}strand.
\newblock Wiener type algebras of pseudodifferential operators.
\newblock In {\em S\'eminaire sur les \'Equations aux D\'eriv\'ees Partielles,
  1994--1995}, pages Exp.\ No.\ IV, 21. \'Ecole Polytech., Palaiseau, 1995.

\bibitem{strohmer04}
T.~Strohmer.
\newblock On the role of the {H}eisenberg group in wireless communication.
\newblock {\em Manuscript}, 2004.

\bibitem{tachizawa94}
K.~Tachizawa.
\newblock The boundedness of pseudodifferential operators on modulation spaces.
\newblock {\em Math. Nachr.}, 168:263--277, 1994.

\bibitem{Toft01}
J.~Toft.
\newblock Subalgebras to a {W}iener type algebra of pseudo-differential
  operators.
\newblock {\em Ann. Inst. Fourier (Grenoble)}, 51(5):1347--1383, 2001.

\bibitem{Toft02a}
J.~Toft.
\newblock Continuity properties in non-commutative convolution algebras, with
  applications in pseudo-differential calculus.
\newblock {\em Bull. Sci. Math.}, 126(2):115--142, 2002.

\bibitem{Toft04}
J.~Toft.
\newblock Continuity properties for modulation spaces, with applications to
  pseudo-differential calculus. {I}.
\newblock {\em J. Funct. Anal.}, 207(2):399--429, 2004.

\bibitem{walnut93}
D.~F. Walnut.
\newblock Lattice size estimates for {G}abor decompositions.
\newblock {\em Monatsh. Math.}, 115(3):245--256, 1993.

\end{thebibliography}

\end{document}